\documentclass[reqno,11pt,a4paper]{amsart}

\usepackage[utf8]{inputenc}
\usepackage[T1]{fontenc}
\usepackage{ae}
\usepackage{units}
\usepackage{icomma}
\usepackage{color}
\usepackage{graphicx}
\usepackage{cancel}
\usepackage{fullpage}
\usepackage{upgreek}
\usepackage{algorithm}
\usepackage{algorithmic}
\usepackage{type1cm}
\usepackage{eso-pic}
\usepackage[breaklinks,pdfpagelabels=false]{hyperref}
\usepackage{verbatim}
\usepackage{tikz}
\usetikzlibrary{arrows}

\usepackage{caption}
\newcommand{\bm}[1]{\mbox{\boldmath $ {#1} $}}

\renewcommand{\epsilon}{\varepsilon}

\newtheorem{theorem}{Theorem}[section]

\newtheorem{remark}[theorem]{Remark}
\theoremstyle{definition}

\numberwithin{equation}{section}
\numberwithin{theorem}{section}

\makeatletter

\def\@settitle{\begin{center}\baselineskip14\p@\relax\bfseries\Large\@title\end{center}}

\makeatother

\begin{document}

\title[Non-freezing of the proximity digraph and non-pseudostable convergence in heterogeneous Hegselmann--Krause models: Counterexamples to conjectures and results of Mirtabatabaei and Bullo] 
      {Non-freezing of the proximity digraph and non-pseudostable convergence in heterogeneous Hegselmann--Krause models: Counterexamples to conjectures and results of Mirtabatabaei and Bullo}



      \date{\today}


      \author{Peter Hegarty, Damiano Ognissanti and Edvin Wedin} \address{Department of Mathematical Sciences, Chalmers University of Technology and University of Gothenburg, 41296 Gothenburg, Sweden} 
\email{hegarty@chalmers.se, damiano@chalmers.se, wedin@fripost.org}







\begin{abstract}
  In 2012, Mirtabatabaei and Bullo studied heterogeneous
  Hegselmann--Krause models of opinion dynamics, where different agents may have
  different confidence or influence bounds. They conjectured that opinion
  vectors always converge, a fundamental problem that remains open. In support
  of this main conjecture they gave some partial results and made some
  auxiliary conjectures. Here we give counterexamples to two of these auxiliary
  conjectures, and even to one theorem in their paper. We discuss to what extent the dynamics proposed in these disproved conjectures and false theorem can nevertheless be expected to hold for randomly chosen initial opinion vectors.
  In the last section (Section \ref{sect:AI}), we describe the role of AI in the production
  of the results in this paper. 
\end{abstract}


\maketitle

\setcounter{equation}{0}

\setcounter{equation}{0}
\setcounter{section}{0}

$\;$ \\
$\;$ \par

\section{{\bf Introduction}}\label{sect:intro}

The Hegselmann--Krause model \cite{HK} is perhaps the
single most studied model in the field of opinion dynamics. Lorenz \cite{L} introduced a \emph{heterogeneous} version of the model, whose
basic idea is that different agents may be more or less open-minded.  Mirtabatabaei and
Bullo \cite{MB} extended this idea to take account of the fact that different agents may be more or less influential. 

Specifically, given a finite set of $n$ agents, their
initial opinions are modelled as a vector of real numbers
$\bm{x}(0) = (x_1(0),\,\dots,\,x_n(0)) \in \mathbb{R}^n$, with the convention
being that $x_1(0) \leq \dots \leq x_n(0)$. Each agent
$i \in \{1,\,\dots,\,n\}$ possesses a fixed positive real number $r_i$, termed
its confidence or influence bound. Agents update their opinions synchronously
in discrete time according to the formula
\begin{equation}\label{eq:dynamics}
  x_i (t+1) = \frac{1}{|\mathcal{N}_i (t)|} \sum_{j \in \mathcal{N}_i(t)} x_j (t), 
\end{equation}
where the so-called \emph{neighbourhood} of agent $i$ at time $t$ is defined as
follows:
\\
\\
For the synchronized bounded confidence (SBC) model (due to Lorenz):
\begin{equation*}
  \mathcal{N}_{i}(t) =\{j: |x_j (t) - x_i(t)| \leq r_i\}.
\end{equation*}
For the synchronized bounded influence (SBI) model (due to Mirtabatabaei and Bullo):
\begin{equation*}
  \mathcal{N}_{i}(t) =\{j: |x_j (t) - x_i(t)| \leq r_j\}.
\end{equation*}
In words, in the SBC model each agent updates its opinion to the average of
those in which it has confidence at time $t$, where having confidence in an
opinion means that it lies within the agent's confidence bound of its own
current opinion. In particular, an agent always retains confidence in its own
opinion. In the SBI model, each agent updates to the average of those who
influence it at time $t$, where the agent is influenced by any other whose
current opinion differs from its own by at most the influence bound of the
other agent. Thus, here also an agent always remains influenced by itself.

Mirtabatabaei and Bullo conjectured (Conjecture 2.1 in \cite{MB}) that, for any
initial opinion vector $\bm{x}(0) \in \mathbb{R}^n$, the vectors $\bm{x}(t)$
converge (in the usual Euclidean topology, thus pointwise) as
$t \rightarrow \infty$. This conjecture remains open and is definitely one of
the outstanding open problems concerning Hegselmann--Krause dynamics. The
original model proposed by Hegselmann and Krause \cite{HK} is the special,
so-called \emph{homogeneous}, case when all $r_i$ are equal, in which case the
SBC and SBI models coincide. This special case has important properties which
do not carry over to the general case, in particular:
\par (1) in the homogeneous case, there is always a finite time $\tau$ at which
the opinion vector \emph{freezes}, i.e.: $\bm{x}(t) = \bm{x}(\tau)$ for all
$t \geq \tau$,
\par (2) the opinion vector remains ordered, i.e.: if $x_i(0) \leq x_j(0)$,
then $x_i(t) \leq x_j (t)$ for all $t$.
\\
Failure of these properties in the general SBC or SBI models is well-known and
it seems to be the basic reason why proving (or disproving) convergence in
general is so difficult.

The paper \cite{MB} can be considered as a kind of roadmap towards attacking
the core problem of convergence. It presents a number of auxiliary conjectures
and some partial results. We now describe two of those conjectures and one of
the main results.

Firstly, for any Hegselmann--Krause system on $n$ agents, its
\emph{proximity digraph} $G$ is the digraph whose nodes are the $n$ agents,
with an arc $(i,\,j)$ if and only if $j \in \mathcal{N}_i$. For
a system evolving according to (\ref{eq:dynamics}), $G_t$ represents the
proximity digraph at time $t$. Mirtabatabaei and Bullo conjectured the
following:
\\
\\
{\bf Conjecture 2.2 in \cite{MB}.} \emph{For any SBC or SBI system evolving
  according to (\ref{eq:dynamics}), there exists a finite time $\tau$ such that
  the proximity digraph $G_t$ is constant for all $t \geq \tau$.}
\\
\par The idea here seems to be to break a proof of Conjecture 2.1 into two
steps, first to prove Conjecture 2.2 and then deduce Conjecture 2.1 from it.
One of the main results in the paper, Theorem 6.4, claims to accomplish the
second step of this plan. Part (iv) of that theorem asserts that, if the
proximity digraph of an SBC or SBI system remains constant after finite time, then either the system freezes in finite time or it converges in a particular manner which, following \cite{U}, they term \emph{pseudostable}.{\footnote{In fact, a literal reading of Theorem 6.4(iv) in \cite{MB} excludes the case of freezing in finite time. We assume this is just a misprint/oversight and that the authors meant to include that option.}}
Pseudostable behaviour means the agents can be partitioned into two non-empty subsets, of freezing or monotonically converging agents. More precisely, there is some finite $\tau$ such that, for every $i \in \{1,\,\dots,\,n\}$, either 
\par (a) $x_i(t) = x_i(t+1)$  for all $t \geq \tau$, or 
\par (b.i) $x_i(t) < x_i (t+1) < x_i (\infty)$ for all $t \geq \tau$, or
\par (b.ii) $x_i (t) > x_i(t+1) > x_i(\infty)$ for all $t \geq \tau$.
\\
Option (a) says that an agent freezes in finite time, whereas options (b) say that it monotonically approaches its limiting position from one side without ever reaching it. Pseudostability further requires there to be at least one agent of each type. 

Finally, Conjecture 2.3 of \cite{MB} asserts that every SBC or SBI system either freezes in finite time or exhibits pseudostable behaviour. This would follow from Conjecture 2.2 and
Theorem 6.4(iv). We note that, in 2017, Su et al \cite{SY} proved that, in any SBC system, the extremal opinions always become frozen in finite time. Their proof easily carries over to SBI systems too. Hence, what has remained open in Conjecture 2.3 is whether every agent that doesn't freeze in finite time must satisfy option (b.i) or (b.ii).
\\
\par Below we will provide counterexamples to all the above assertions. The
rest of the paper is divided into six short sections. Sections \ref{sect:SBC2}
and \ref{sect:SBI2} present counterexamples to Conjectures 2.2 and 2.3 for the
SBC and SBI models respectively. Sections \ref{sect:SBC6} and \ref{sect:SBI6}
provide counterexamples to parts (iii) and (iv) of Theorem 6.4 for the
respective models (part (iii) is rather technical and we refer the reader to
\cite{MB} for its precise statement). Section \ref{sect:mod} discusses the
remaining plausibility of the approach in \cite{MB} for most initial
configurations. For while our examples
demonstrate that the roadmap outlined in \cite{MB} is fundamentally flawed, we
will explain why we think that, at least for the part concerning Conjecture 2.2, it may still be plausible for a ``generic'' initial configuration of opinions, something which is already
supported by what can be recovered from \cite{MB} itself. Finally,
Section \ref{sect:AI} discusses how AI was used in producing the results of
this paper.

\vspace{3mm}

\section{{\bf Counterexample to Conjectures 2.2 and 2.3 for the SBC model}}\label{sect:SBC2}

We take $n=7$, initial opinion vector
\begin{equation*}
  \bm{x}(0) = (x_1(0),\,\dots,\,x_7(0)) = (0,\,38,\,69,\,84,\,99,\,130,\,168)
\end{equation*}
and confidence radii
\begin{equation*}
  \bm{r} = (r_1,\,\dots,\,r_7) = (18,\,42,\,48,\,12,\,48,\,42,\,18).
\end{equation*}
We claim that, for all $t \geq 0$,
\begin{equation}\label{eq:F1}
  \bm{x}(t) = (0,\,36,\,72,\,84,\,96,\,132,\,168) + \left( - \frac{1}{6} \right)^t \cdot (0,\,2,\,-3,\,0,\,3,\,-2,\,0).
\end{equation}
This is directly verified to be the case for $t=0$, so we can proceed by
induction on $t$. Note that the initial configuration is symmetric about
position $x=84$. Since the confidence radii are also symmetric about the middle
agent $\# 4$, it is immediate that $\bm{x}(t)$ will remain symmetric about
$x = 84$ for all $t$. Assume (\ref{eq:F1}) holds at some time $t \geq 0$. It
thus suffices to verify that it still holds at time $t+1$ for $x_i(t+1)$ and
$1 \leq i \leq 3$. We take these three agents one at a time.
\\
\\
$i = 1$: Clearly $x_2(t) > 18$, so agent 1 only listens to itself and
$x_1 (t+1) = x_1(t) = 0$, as desired.
\\
\\
$i = 2$: It is easy to check that agent 2 will listen to 1, 2 and 3 at time
$t$ if (\ref{eq:F1}) holds. Thus
\begin{eqnarray*}
  x_2 (t+1) = \frac{1}{3} \sum_{j=1}^{3} x_j (t) = \frac{1}{3} (0+36+72) + \frac{1}{3} \left( -\frac{1}{6} \right)^t (0+2-3) = 36 + 2 \cdot \left( - \frac{1}{6} \right)^{t+1},
\end{eqnarray*}
as desired.
\\
\\
$i=3$: It is easy to check that agent 3 will listen to 2, 3, 4 and 5 if
(\ref{eq:F1}) holds. Thus
\begin{eqnarray*}
  x_3 (t+1) = \frac{1}{4} \sum_{j=2}^{5} x_j (t) = \frac{1}{4} (36+72+84+96) + \frac{1}{4} \left( -\frac{1}{6} \right)^t (2-3+0+3) = 72 - 3 \cdot \left( - \frac{1}{6} \right)^{t+1},
\end{eqnarray*}
as desired. This completes the induction. 
\\
\\
If $t$ is even, then (\ref{eq:F1}) implies that $x_3 (t) < 72$ and
$x_5 (t) > 96$. In this case, agent 4 listens only to itself and thus
$x_4 (t+1) = x_4 (t) = 84$, as desired. If $t$ is odd, then $x_3(t) > 72$ and
$x_5(t) < 96$, so agent 4 listens to both 3 and 5. Formula (\ref{eq:F1})
implies that $\bm{x}(t)$ converges to $(0,\,36,\,72,\,84,\,96,\,132,\,168)$ as
$t \rightarrow \infty$. But the proximity digraph never freezes, since agent 4
listens to both 3 and 5
at time $t$ if and only if $t$ is odd. Thus, we have a counterexample to
Conjecture 2.2 in \cite{MB}. Moreover, according to (\ref{eq:F1}), agents
2, 3, 5 and 6 continually hop between opposite sides of their limiting
locations. Thus, we also have a counterexample to Conjecture 2.3 in \cite{MB}.

\begin{remark}\label{rem:wiggle}
  There is obviously some wiggle-room in the choice of the precise numbers in
  the components of the vectors $\bm{x}(0)$ and $\bm{r}$ above. We choose not
  to delve into this matter, either for this or for subsequent examples below. See
  Section \ref{sect:AI} however for some further remarks on ``classification''
  of counterexamples to the assertions in \cite{MB}.
\end{remark}

\vspace{3mm}

\section{{\bf Counterexample to Conjectures 2.2 and 2.3 for the SBI model}}\label{sect:SBI2}

We take $n=7$ again. Denote
\begin{equation*}
  \phi = \frac{1+\sqrt{5}}{2}, \;\;\; \lambda = \frac{1-\sqrt{5}}{8}
\end{equation*}
and observe that
\begin{equation}\label{eq:F2}
  \lambda = \frac{1-\phi}{4} = - \frac{1}{4\phi}.
\end{equation}
Choose as the initial opinion vector
\begin{equation*}
  \bm{x}(0) = (x_1(0),\,\dots,\,x_7(0)) = \left( 0,\,\frac{15}{2},\, 10 - \frac{\phi}{2},\, 11,\, 12 + \frac{\phi}{2},\, \frac{29}{2},\,22 \right)
\end{equation*}
and influence radii
\begin{equation*}
  \bm{r} = (r_1,\,\dots,\,r_7) = (8,\,4,\,4,\,8,\,4,\,4,\,8).
\end{equation*}
We claim that, for all $t \geq 0$,
\begin{equation}\label{eq:F3}
  \bm{x}(t) = (0,\,7,\,10,\,11,\,12,\,15,\,22) + \frac{\lambda^t}{2} \cdot (0,\,1,\,-\phi,\,0,\,\phi,\,-1,\,0).
\end{equation}
This is directly verified to be the case for $t=0$, so assume (\ref{eq:F3})
holds at time $t$. The initial configuration is symmetric about position $x=11$
and the influence radii are symmetric about the middle agent $\# 4$. Thus it
again suffices to verify that (\ref{eq:F3}) still holds at time $t+1$ for
agents 1, 2 and 3. We take them one at a time.
\\
\\
$i = 1$: Clearly $x_2(t) > 6$ and, since $r_2  = 4$, agent 1 will not be
influenced by 2. It's easily checked that it won't be influenced by any other
agent either. Thus $x_1 (t+1) = x_1(t) = 0$, as desired.
\\
\\
$i = 2$: It is easy to check that agent 2 will listen to 1, 2, 3 and 4 at time
$t$ if (\ref{eq:F3}) holds. Thus
\begin{eqnarray*}
  x_2 (t+1) = \frac{1}{4} \sum_{j=1}^{4} x_j (t) = \frac{1}{4} (0+7+10+11) + \frac{1}{4} \cdot \frac{\lambda^t}{2} \cdot (0+1-\phi+0) = 7 + \frac{\lambda^{t+1}}{2},
\end{eqnarray*}
where we have used the first equality in (\ref{eq:F2}).
\\
\\
$i=3$: It is easy to check that agent 3 will listen to 2, 3, 4 and 5 if
(\ref{eq:F3}) holds. Thus
\begin{eqnarray*}
  x_3 (t+1) = \frac{1}{4} \sum_{j=2}^{5} x_j (t) = \frac{1}{4} (7+10+11+12) + \frac{1}{4} \cdot \frac{\lambda^t}{2} \cdot (1-\phi+0+\phi) = 10 - \phi \cdot \frac{\lambda^{t+1}}{2},
\end{eqnarray*}
where we have used the second equality in (\ref{eq:F2}). This completes the
induction. 
\\
\\
Note that $\lambda < 0$. If $t$ is odd, then (\ref{eq:F3}) implies that
$x_2 (t) < 7$ and $x_6 (t) > 15$. In this case, agent 4 listens to 3, 4 and 5.
If $t$ is even, then $x_2(t) > 7$ and $x_6(t) < 15$, so agent 4 will also
listen to 2 and 6. Formula (\ref{eq:F3}) implies that $\bm{x}(t)$ converges to
$(0,\,7,\,10,\,11,\,12,\,15,\,22)$ as $t \rightarrow \infty$. But the
proximity digraph never freezes, since agent 4 listens to both 2 and 6 at
time $t$ if and only if $t$ is even. Moreover, since $\lambda < 0$, agents
2, 3, 5 and 6 continually hop between opposite sides of their limiting
locations. Thus, we again have a counterexample to Conjectures 2.2 and 2.3
in \cite{MB}.

\vspace{3mm}

\section{{\bf Counterexamples to Theorem 6.4 for the SBC model}}\label{sect:SBC6}

We take $n=6$, initial opinion vector
\begin{equation*}
  \bm{x}(0) = (x_1(0),\,\dots,\,x_6(0)) = (0,\,22,\,37,\,103,\,118,\,140)
\end{equation*}
and confidence radii
\begin{equation*}
  \bm{r} = (r_1,\,\dots,\,r_6) = (10,\,70,\,70,\,70,\,70,\,10).
\end{equation*}
Set $\xi = -\frac{1}{6}$. We claim that, for all $t \geq 0$,
\begin{equation}\label{eq:F4}
  \bm{x}(t) = (0,\,20,\,40,\,100,\,120,\,140) + \xi^t \cdot (0,\,2,\,-3,\,3,\,-2,\,0).
\end{equation}
This is directly verified for $t=0$, so we proceed by induction on $t$. The
initial configuration is symmetric about position $x=70$ and the confidence
radii are symmetric about $i=3.5$. So assuming (\ref{eq:F4}) holds at time $t$,
it suffices to verify that it still does so at time $t+1$ for agents
1, 2 and 3. We take them one at a time.
\\
\\
$i = 1$: Clearly $x_2(t) > 10$, so agent 1 only listens to itself and
$x_1 (t+1) = x_1(t) = 0$, as desired.
\\
\\
$i = 2$: It is easy to check that agent 2 will listen to 1, 2 and 3 at time $t$
if (\ref{eq:F4}) holds. Thus
\begin{eqnarray*}
  x_2 (t+1) = \frac{1}{3} \sum_{j=1}^{3} x_j (t) = \frac{1}{3} (0+20+40) + \frac{1}{3} \xi^t (0+2-3) = 20 + 2 \cdot \xi^{t+1},
\end{eqnarray*}
as desired.
\\
\\
$i=3$: It is easy to check that agent 3 will listen to 1, 2, 3 and 4 if
(\ref{eq:F4}) holds. Thus
\begin{eqnarray*}
  x_3 (t+1) = \frac{1}{4} \sum_{j=1}^{4} x_j (t) = \frac{1}{4} (0+20+40+100) + \frac{1}{4} \xi^t (0+2-3+3) = 40 - 3 \cdot \xi^{t+1},
\end{eqnarray*}
as desired. This completes the induction.

Formula (\ref{eq:F4}) is easily checked to imply that the proximity digraph
remains the same for all $t$ and consists, apart from all arcs $(i,\,i)$,
$1 \leq i \leq 6$, of the following arcs (recall that an arc $(i,\,j)$
indicates that agent $i$ takes account of the opinion of $j$):
\begin{equation*}
  (2,\,1), \, (3,\,1), \, (3,\,2), \, (2,\,3), \, (4,\,3), \, (3,\,4), \, (5,\,4), \, (4,\,5), \, (4,\,6), \, (5,\,6).
\end{equation*}
But, since $\xi < 0$, agents 2, 3, 4 and 5 continually hop between opposite
sides of their limiting locations. This contradicts part (iv) of Theorem 6.4.
Moreover, the per-step convergence factor, as defined in Definition 6.1 of
\cite{MB}, is obviously equal to $\xi$. This already obviously contradicts part
(iii)(a) of Theorem 6.4, since $\xi < 0$. Agents 2--5 form an open-minded SCC
with adjacency matrix
\begin{equation}\label{eq:M1}
  \Theta = \left( \begin{array}{cccc} 1/3 & 1/3 & 0 & 0 \\ 1/4 & 1/4 & 1/4 & 0 \\ 0 & 1/4 & 1/4 & 1/4 \\ 0 & 0 & 1/3 & 1/3 \end{array} \right).
\end{equation}
It is easily checked that the eigenvalues of this matrix are
$-\frac{1}{6}, \, \frac{1}{2}, \, \frac{5 \pm \sqrt{13}}{12}$, and hence the
spectral radius is $R = \frac{5 + \sqrt{13}}{12} \approx 0.7171 ...$.
Theorem 6.4(iii)(a) says that the per-step convergence factor converges to $R$,
which is not the case here. 
\\
\par We now turn to part (iii)(b) of Theorem 6.4. The hypothesis there assumes
that the proximity digraph contains at least two distinct open-minded SCCs,
so the example just presented doesn't apply. However, we can modify it to get a
counterexample to (iii)(b) by adding a few ``passive'' agents. Precisely, we
add three more agents to the right of the existing six. Agent 7, called the
\emph{spectator}, starts at $x_7(0) = 210$ and has $r_7 = 100$. Agents 8 and 9,
called the \emph{beacons}, start even further to the right at $x_8(0) = 280$
and $x_9(0) = 300$. They satisfy $r_8 = r_9 = 10$. It is then easy to check, by
induction on $t$, that the following hold:
\\
\\
{\sc Claim 1:} None of agents 1--6 ever listen to any of 7, 8 or 9. In return,
agents 8 and 9 never listen to anyone but themselves (not even each other). But
agent 7 always listens to 8 and 9, as well as to 5 and 6.
\\
\\
{\sc Claim 2:} Hence, $x_1(t),\,\dots,\,x_6(t)$ evolve precisely as in
(\ref{eq:F4}), while $x_8 (t) = 280$ and $x_9 (t) = 300$ for all $t$.
\\
\\
{\sc Claim 3:} $x_7 (t) = 210 + \frac{12}{11} \left( \xi^t - \left( \frac{1}{5} \right)^t \right)$.
\\
\par Claims 1 and 2 are verified essentially by repeating the previous
computations, so let us verify Claim 3. Suppose it and Claims 1 and 2 hold at
time $t$. Then
\begin{eqnarray*}
  x_7 (t+1) = \frac{1}{5} \sum_{j=5}^{9} x_j (t) = \frac{1}{5} (120 + 140 + 210 + 280 + 300) + \frac{1}{5} \left[ - 2 \cdot \xi^t + \frac{12}{11} \cdot \xi^t - \frac{12}{11} \cdot \left( \frac{1}{5} \right)^t \right] \\ = 210 - \frac{12}{11} \cdot \left( \frac{1}{5} \right)^{t+1} - \frac{2}{11} \cdot \xi^t = 210 + \frac{12}{11} \left( \xi^{t+1} - \left( \frac{1}{5} \right)^{t+1} \right), \;\; {\hbox{since $\xi = -\frac{1}{6}$}}.
\end{eqnarray*}

Now Claim 1 implies that, in the notation and terminology of \cite{MB}, there
are two open-minded SCCs in this 9-agent configuration, that containing agents
2--5 as before plus that containing the single agent 7. The former is the
successor of the latter, since 7 listens to 5, but nobody
listens to 7. Thus, $S_k = \{7\}$ and $S_m = \{2,\,3,\,4,\,5\}$. Since nobody
in $S_m$ listens to 7, 8 or 9, we have $\Theta_m = \Theta$, the matrix in
(\ref{eq:M1}). Since 7 listens to five agents,
$\Theta_k = \left( \frac{1}{5} \right)$. Thus $\Theta_k$ has spectral radius
$1/5$, strictly less than the previously noted spectral radius
$\frac{5 + \sqrt{13}}{12}$ of $\Theta_m$. The hypotheses of Theorem 6.4(iii)(b)
are therefore satisfied. As previously noted, since $\xi < 0$, agents 2--5
continually hop between opposite sides of their limiting locations. Moreover,
at even $t$, 2 and 4 are to the right of their limiting locations whereas 3 and
5 are to the left. The situation is reversed at odd $t$. The only way this
could be consistent with the conclusion of Theorem 6.4(iii)(b) is if agent 7
froze at $x_7 = 210$ from some time onwards. But Claim 3 implies agent 7 is
strictly to the left of 210 for all $t \geq 1$. 

\vspace{3mm}

\section{{\bf Counterexamples to Theorem 6.4 for the SBI model}}\label{sect:SBI6}

We take $n=7$. Let $\phi$ and $\lambda$ be as in Section \ref{sect:SBI2}. We
choose as the initial opinion vector
\begin{equation*}
  \bm{x}(0) = (x_1(0),\,\dots,\,x_7(0)) = (0,\,71,\,100-\phi,\,110,\,120+\phi,\,149,\,220)
\end{equation*}
and influence radii
\begin{equation*}
  \bm{r} = (r_1,\,\dots,\,r_7) = (85,\,35,\,35,\,75,\,35,\,35,\,85).
\end{equation*}
We claim that, for all $t \geq 0$,
\begin{equation}\label{eq:F5}
  \bm{x}(t) = (0,\,70,\,100,\,110,\,120,\,150,\,220) + \lambda^t \cdot (0,\,1,\,-\phi,\,0,\,\phi,\,-1,\,0).
\end{equation}
As usual, this is first verified for $t=0$ and we proceed by induction on $t$.
The initial configuration is symmetric about position $x=110$ and the influence
radii are symmetric about $i=4$. So assuming (\ref{eq:F5}) holds at time $t$,
it suffices to verify that it still does so at time $t+1$ for agents
1, 2 and 3. We take them one at a time.
\\
\\
$i = 1$: Clearly $x_2(t) > r_2 = 35$, so agent 1 is not influenced by 2. It is
easily checked that it's not influenced by any other agent but itself, and so
$x_1 (t+1) = x_1(t) = 0$, as desired.
\\
\\
$i = 2$: It is easy to check that agent 2 will be influenced by 1, 2, 3 and 4
at time $t$ if (\ref{eq:F5}) holds. Thus
\begin{eqnarray*}
  x_2 (t+1) = \frac{1}{4} \sum_{j=1}^{4} x_j (t) = \frac{1}{4} (0+70+100+110) + \frac{1}{4} \lambda^t (0+1-\phi) = 70 + \lambda^{t+1},
\end{eqnarray*}
where we have used the first equality in (\ref{eq:F2}).
\\
\\
$i=3$: It is easy to check that agent 3 will be influenced by 2, 3, 4 and 5 if
(\ref{eq:F5}) holds. Thus
\begin{eqnarray*}
  x_3 (t+1) = \frac{1}{4} \sum_{j=2}^{5} x_j (t) = \frac{1}{4} (70+100+110+120) + \frac{1}{4} \lambda^t (1-\phi+0+\phi) = 100 - \phi \cdot \lambda^{t+1},
\end{eqnarray*}
where we have used the second equality in (\ref{eq:F2}). This completes the
induction.

Formula (\ref{eq:F5}) is easily checked to imply that the proximity digraph
remains the same for all $t$ and consists, apart from all arcs $(i,\,i)$,
$1 \leq i \leq 7$, of the following arcs (as for the SBC model, an arc
$(i,\,j)$ indicates that agent $i$ takes account of the opinion of $j$):
\begin{equation*}
  (2,\,1), \, (3,\,2), \, (2,\,3), \, (4,\,3), \, (5,\,3), \, (2,\,4), \, (3,\,4), \, (5,\,4), \, (6,\,4), \, (3,\,5), \, (4,\,5), \, (6,\,5), \, (5,\,6), \, (6,\,7).
\end{equation*}
Since $\lambda < 0$, agents 2, 3, 5 and 6 continually hop between opposite
sides of their limiting locations, which contradicts part (iv) of Theorem 6.4.
The per-step convergence factor is obviously equal to $\lambda$, immediately
contradicting part (iii)(a) of Theorem 6.4, since $\lambda < 0$. Agents 2--6
form an open-minded SCC with adjacency matrix
\begin{equation}\label{eq:M2}
  \Theta = \left( \begin{array}{ccccc} 1/4 & 1/4 & 1/4 & 0 & 0 \\ 1/4 & 1/4 & 1/4 & 1/4 & 0 \\ 0 & 1/3 & 1/3 & 1/3 & 0 \\ 0 & 1/4 & 1/4 & 1/4 & 1/4 \\ 0 & 0 & 1/4 & 1/4 & 1/4 \end{array} \right).
\end{equation}
It is easily checked that the largest eigenvalue of this matrix is
$R \approx 0.9533$. Theorem 6.4(iii)(a) says that the per-step convergence
factor converges to $R$, which is not the case here.
\\
\par We finally turn to part (iii)(b) of Theorem 6.4, where once again the idea
will be to suitably add some ``passive'' agents to the previous configuration so
as to create two different open-minded SCCs and thus satisfy the hypotheses of
(iii)(b). This time it suffices to add a single spectator, we don't need
beacons -- see Remark \ref{rem:beacons} below. The spectator is initially placed
at $x = 60$ and has a tiny influence radius, say $r=1$. Since the convention is
that the components of the initial opinion vector are written in increasing
order, the spectator now becomes agent 2 and what were previously agents 2--7
now become 3--8. The following are then verified by induction on $t$, as usual.
\\
\\
{\sc Claim 1:} Agent 2 never influences anyone else, but 2 is always under the
influence of 1, 3 and 5.
\\
\\
{\sc Claim 2:} The agents other than 2 evolve as in (\ref{eq:F5}).
\\
\\
{\sc Claim 3:} $x_2 (t) = 60 + \frac{1}{\phi} \left[ \left( \frac{1}{4} \right)^t - \lambda^t \right]$.
\\
\par Let us complete the induction step for Claim 3. Suppose Claims 1, 2 and 3
all hold at time $t$. Then
\begin{eqnarray*}
  x_2 (t+1) = \frac{1}{4} (x_1(t) + x_2(t) + x_3(t) + x_5(t)) = \frac{1}{4} (0+60+70+110) + \frac{1}{4} \left[ \frac{1}{\phi} \left( \left( \frac{1}{4} \right)^t - \lambda^t \right) + \lambda^t \right] \\ = 60 + \frac{1}{\phi} \cdot \left( \frac{1}{4} \right)^{t+1} - \frac{\lambda^t}{4} \left( \frac{1}{\phi} - 1 \right) = 60 + \frac{1}{\phi} \left[ \left( \frac{1}{4} \right)^{t+1} - \lambda^{t+1} \right],
\end{eqnarray*}
as desired, where we have used the first equality in (\ref{eq:F2}) again.
\\
\par The previous 5-agent open-minded SCC remains, now formed by agents 3--7.
Since none of them ever listen to the spectator 2, the adjacency matrix of this
SCC is, as before, the matrix $\Theta$ in (\ref{eq:M2}), and so
$R \approx 0.9533$. Agent 2 forms an open-minded SCC of its own and, since it
always listens to three other agents, two of whom are in the other open-minded
SCC, its spectral radius is $1/4$. So the hypotheses of Theorem 6.4(iii)(b) are
satisfied and agent 2 is supposed to follow its leader. But, since
$\lambda < 0$, we run into exactly the same problem as in Section
\ref{sect:SBC6}. The only way for the conclusion of Theorem 6.4(iii)(b) to be
consistent with the oscillations of agents 3, 4, 6 and 7 would be for agent 2 to become frozen at $x_2 = 60$ from some time onwards. But Claim 3 implies that agent 2 is strictly to the right of 60 at all $t \geq 1$. 

\begin{remark}\label{rem:beacons}
  In the influence model, we could place the spectator in amongst existing
  agents and ensure it didn't affect their behaviour just by giving it a tiny
  influence radius. In the confidence model this doesn't work because if any
  arc of the pre-existing proximity digraph crosses the spectator's location
  then the source of the arc must necessarily have confidence in the spectator
  also. Hence, in the confidence model we needed to place the spectator far to
  the right, beyond the confidence bounds of all existing agents. The beacons
  then needed
  to be added to ``stabilize'' the spectator, i.e.: to ensure that it remained
  close to its initial position. Two beacons were needed because we wanted the
  spectator to interact with at least one agent in the existing open-minded
  SCC. 
\end{remark}

\vspace{3mm}

\section{{\bf Modifications of the assertions in Mirtabatabaei--Bullo}}\label{sect:mod}

We start with the example in Section \ref{sect:SBC2}. The configuration
converges to
\begin{equation*}
  \bm{x}_{\infty} = (0,\,36,\,72,\,84,\,96,\,132,\,168).
\end{equation*}
In this limiting configuration, agents 3 and 5 are at distance exactly
$r_4 = 12$ from agent 4. In the notation and terminology of Definition 4.1 of
\cite{MB}, the 4th component of the equitopology distance of $\bm{x}_{\infty}$
is zero, i.e.: $\epsilon_4(\bm{x}_{\infty}) = 0$.

Similarly, in the example of Section \ref{sect:SBI2}, agent 4 sits at
$x_4 = 11$ while agents 2 and 6 converge to $x_2 = 7$ and $x_6 = 15$
respectively. Since $r_2 = r_6 = 4$, the limiting configuration satisfies
$\epsilon_2 = \epsilon_6 = 0$.

Now Lemma 4.8 of \cite{MB} asserts that, in both the SBC and SBI models, if an
initial configuration converges and every component of the equitopology
distance of the limiting configuration is non-zero, then Conjecture 2.2 holds. It is thus eminently plausible that Conjecture
2.2 holds for a \emph{generic} initial configuration (that is, holds with
probability one when the initial configuration is chosen uniformly at random
from a suitable probability space). 
Of course, Lemma 4.8 assumes
convergence in the first place, which suggests that, while the ``roadmap''
sketched in \cite{MB} for understanding the heterogeneous SBC and SBI models in
full generality fails, it may still be valid for generic initial configurations
of opinions.
\\
\par Turning to Theorem 6.4 of \cite{MB}, consider first the example of Section
\ref{sect:SBC6}. Agents 2--5 form an open-minded SCC whose opinion vector
converges to $(20,\,40,\,100,\,120)$. By (\ref{eq:F4}), the offset from this
limiting state at $t=0$ is $\bm{v} = (2,\,-3,\,3,\,-2)$.

The adjacency matrix $\Theta$ of the open-minded SCC has largest eigenvalue
$R = \frac{5+\sqrt{13}}{12}$ and its left-eigenspace is spanned by a
vector of the form $\bm{w} = (\alpha,\,\beta,\,\beta,\,\alpha)$. It follows
that $\bm{v} \cdot \bm{w} = 0$. So the initial offset $\bm{v}$ of the
open-minded SCC is orthogonal to the principal left-eigenvector of the
adjacency matrix of that SCC. This is the reason why the proof of
Theorem 6.4(iii)(a) fails for this example. 

A similar analysis applies to the other examples in Sections \ref{sect:SBC6}
and \ref{sect:SBI6}, and we omit the details. Note, however, that the hypothesis of part (iii)(b) of Theorem 6.4 also assumes there is a unique leader SCC in $\mathcal{M}(S_k (y))$ whereas, in general, there could be multiple SCCs with equal spectral radii. Indeed, Definition 6.3 of \cite{MB} does not uniquely determine a ``leader" SCC in that case. On the other hand, part (iv) does not depend on (iii)(b), only on (iii)(a). 
The upshot is that Theorem 6.4(iv) may still hold for generic
initial configurations. Finally, the same can be said of Conjecture 2.3, which as previously noted would follow from Conjecture 2.2 and Theorem 6.4(iv). Section 7 of \cite{MB} makes reference to
simulations which support this hypothesis for Conjectures 2.2 and 2.3. Once again, the
takeaway is that the ``roadmap'' of \cite{MB} may still be useful for
understanding most, if not all, possible configurations of opinions. On the
other hand, it suggests that, if the main Conjecture 2.1 is true, proving it
will require an approach not tailored to generic configurations. 

\vspace{3mm}

\section{{\bf AI statement}}\label{sect:AI}

All of the examples above were found by \texttt{Claude Opus 5.5}, upon
prompting by the authors. We asked it whether the proximity digraph must
freeze in finite time, assuming that the opinion vector converges, but without
making any direct reference to the paper \cite{MB}. It returned with the
counterexample in Section \ref{sect:SBC2} and further prompting led to the
other counterexamples above.

The text of this paper is entirely the work of the authors, though we must also
give \texttt{Claude} credit for coming up with the terms ``spectator'' and
``beacon''. In particular, we have verified by hand the correctness of all the
examples presented above.

Our interactions with \texttt{Claude} on this topic have so far also generated
a lot more technical data, for example concerning the minimal $n$ for which
counterexamples exist in the various settings of Sections
\ref{sect:SBC2}--\ref{sect:SBI6}, plus attempts at classification of all possible
counterexamples for such minimal $n$. We do not include such material here, for
two reasons. Firstly, a lot of it is technical and reliant even on large
computations which \texttt{Claude} performed in the background. These have not
been verified in full by us, the human authors. Secondly, and just as
importantly, we feel such technical material would distract from the simple
elegance of the examples presented in Sections \ref{sect:SBC2}--\ref{sect:SBI6},
which provide counterexamples to well-known, longstanding conjectures, as well
as a related published theorem, and crystallize how the material presented in a
highly cited paper needs to be revised.

\vspace{3mm}

\section*{{\bf Note added post-publication}}

A couple of weeks after the initial version of our preprint appeared on \texttt{arxiv.org}, Etesami and Touri \cite{ET} announced a counterexample to Conjecture 2.1 of \cite{MB}. Their counterexample is for the SBC model and has 6 agents. It follows from their construction that Conjectures 2.1, 2.2 and 2.3 of \cite{MB} all fail for an open set of initial data, and hence the approach of \cite{MB} does not, in fact, hold even for generic initial configurations, at least in the SBC model. At the time of writing, no counterexample to Conjecture 2.1 for the SBI model has yet been publicized. Note that, strictly speaking, the result in \cite{ET} doesn't supersede ours, since our examples show that the various conjectures and results in \cite{MB} may fail even for a HK-system that does converge.  


\end{document}